\documentclass[12pt]{amsart}

\usepackage{amssymb}
\usepackage{epsf}

\newcommand{\Span}[1]{\langle #1 \rangle}

\newcommand{\QED}{
\setlength{\unitlength}{1.0pt}%
\begin{picture}(10,7.5)
\put(0,-2.5){\rule{2.5pt}{2.5pt}}
\put(0,0){\rule{5pt}{2.5pt}}
\put(0,2.5){\rule{10pt}{2.5pt}}
\end{picture}\vspace{10pt}}
 
\newtheorem{thm}{Theorem}[section]
\newtheorem{lemma}[thm]{Lemma}
\newtheorem{cor}[thm]{Corollary}

\newtheorem{ex}[thm]{Example}

\newtheorem{rem}[thm]{Remark}

\begin{document}

\title[Pieri-type formula for maximal isotropic Grassmannians]
{Pieri-type formulas for maximal isotropic Grassmannians via triple
intersections} 

\author{Frank Sottile}
\address{Department of Mathematics\\
        University of Toronto\\
        100 St.~George Street\\
	Toronto, Ontario  M5S 3G3\\
	Canada}
\email{sottile@math.toronto.edu}
\date{\today}
\thanks{Research supported in part by NSF grant
    DMS-90-22140 and NSERC grant OPG0170279}

\begin{abstract}
We give an elementary proof of the Pieri-type formula in the cohomology of a
Grassmannian of maximal isotropic subspaces of an odd 
orthogonal or symplectic vector space.
This proof proceeds by explicitly computing a triple intersection of
Schubert varieties.
The decisive step is an exact description of the intersection of two
Schubert varieties, from which the multiplicities (which are powers of 2) in
the Pieri-type formula are immediately obvious.
\end{abstract}

\maketitle

\section*{Introduction}
The goal of this paper is to give an elementary geometric proof of 
Pieri-type formulas in the cohomology of
Grassmannians of maximal isotropic subspaces of odd orthogonal or symplectic
vector spaces.
For this, we explicitly compute a triple intersection of
Schubert varieties, where one is a special Schubert variety.
Previously, Sert\"oz~\cite{Sertoz} had studied such triple intersections in
orthogonal Grassmannians, but was unable to determine
the intersection multiplicities and obtain a formula.

These multiplicities are either 0 or powers of 2.
Our proof explains them as the intersection multiplicity of
a linear subspace (defining the special Schubert variety) 
with a collection of quadrics and linear subspaces (determined by
the other two Schubert varieties).
This is similar to triple intersection proofs of the classical Pieri
formula ({\em cf.}~\cite{Hodge_intersection}\cite[p. 203]{Griffiths_Harris}%
\cite[\S 9.4]{Fulton_tableaux}) where the multiplicities (0 or 1) count the
number of points in the intersection of linear subspaces. 
A proof of the Pieri-type formula for classical flag
varieties~\cite{Sottile_pieri_schubert} was based upon those ideas.
Similarly, the ideas here  provide a basis for a proof of
Pieri-type formulas in the cohomology of symplectic  flag
varieties~\cite{Bergeron_Sottile_Lagrangian_Pieri}.

These Pieri-type formulas are due to  Hiller and Boe~\cite{Hiller_Boe},
whose proof used the Chevalley formula~\cite{Chevalley91}.
Another proof, using the Leibnitz formula for symplectic and orthogonal
divided differences, was given by Pragacz and
Ratajski~\cite{Pragacz_Ratajski_Operator}. 
These formulas also arise in the theory of projective representations of
symmetric groups~\cite{Schur,Hoffman_Humphreys} as product formulas for
Schur $P$- and $Q$-functions, and were first proven in this context by
Morris~\cite{Morris}. 
The connection of Schur $P$- and $Q$-functions to geometry was established
by Pragacz~\cite{Pragacz_88} (see also~\cite{Pragacz_S-Q}).
\smallskip

In Section 1, we give the basic definitions and state the Pieri-type
formulas in both the orthogonal and symplectic cases, and conclude with an
outline of the proof in the orthogonal case.
Since there is little difference between the proofs in each case, we only
do the orthogonal case in full.
In Section 2, we describe the intersection of two
Schubert varieties, which we use in Section 3 to complete the proof.

\section{The Grassmannian of maximal isotropic subspaces} 

Let $V$ be a ($2n+1$)-dimensional complex vector space equipped with a
non-degenerate symmetric bilinear form $\beta$
and $W$ a $2n$-dimensional complex vector space equipped with a
non-degenerate alternating bilinear form, also denoted
$\beta$.
A subspace $H$ of $V$ or of $W$ is {\em isotropic} if the restriction of 
$\beta$ to $H$ is identically zero.
Isotropic subspaces have dimension at most $n$.
The {\em Grassmannian of maximal isotropic subspaces} $B_n$ or $B(V)$
(respectively $C_n$ or $C(W)$) is the collection of all isotropic
$n$-dimensional subspaces of $V$ (respectively of $W$).
The group 
$\mbox{\em So}_{2n+1}{\mathbb C} = {\rm Aut}(V,\beta)$ acts
transitively on $B_n$ with the stabilizer $P_0$ of a point a maximal
parabolic subgroup associated to the short root, hence
$B_n = \mbox{\em So}_{2n+1}{\mathbb C}/P_0$.
Similarly, $C_n= \mbox{\em Sp}_{2n}{\mathbb C}/P_0$, where $P_0$ is a maximal
parabolic associated to the long root.

Both $B_n$ and $C_n$ are smooth complex manifolds of dimension
$\binom{n+1}{2}$.
While they are not isomorphic if $n>1$, they have identical
Schubert decompositions and Bruhat orders.
Another, deeper, connection is discussed in Remark~\ref{rem:similarities}.
We describe the Schubert decomposition.
For an integer $j$, let $\overline{\jmath}$ denote $-j$.
Choose bases $\{e_{\overline{n}},\ldots,e_n\}$ of $V$ and
$\{f_{\overline{n}},\ldots,f_n\}$ of $W$ for which:
$$
\beta(e_i,e_j)\ =\ \left\{
\begin{array}{ll} 1&\mbox{ \ if } i=\overline{\jmath}\\
                  0&\mbox{ \ otherwise}\end{array}\right.
\quad\mbox{and}\quad
\beta(f_i,f_j)\ =\ \left\{
\begin{array}{ll} j/|j|&\mbox{ \ if } i=\overline{\jmath}\\
                       0            &\mbox{ \ otherwise}\end{array}\right.
.
$$
For example, $\beta(e_1,e_0)=\beta(f_{\overline{2}},f_1)=0$ and 
$\beta(e_0,e_0)=
\beta(f_{\overline{1}},f_1)=-\beta(f_1,f_{\overline{1}})= 1$.
Schubert varieties are determined by  sequences
$$
\lambda:\ n\geq \lambda_1>\lambda_2>\cdots>\lambda_n\geq \overline{n}
$$
whose set of absolute values 
$\{|\lambda_1|,\ldots,|\lambda_n|\}$ equals
$\{1,2,\ldots,n\}$.
Let ${\mathbb{SY}}_n$ denote this set of sequences.
The Schubert variety $X_\lambda$ of $B_n$ is
$$
\{H\in B_n\mid\dim 
H\cap\Span{e_{\lambda_j},\ldots,e_n}\geq j
\mbox{\ for } 1\leq j\leq n\}
$$
and the Schubert variety $Y_\lambda$ of $C_n$
$$
\{H\in C_n\mid\dim 
H\cap\Span{f_{\lambda_j},\ldots,f_n}\geq j
\mbox{\ for } 1\leq j\leq n\}.
$$
Both $X_\lambda$ and $Y_\lambda$ have codimension 
 $|\lambda|:=\lambda_1+\cdots+\lambda_k$, where 
$\lambda_k>0>\lambda_{k+1}$.
Given $\lambda,\mu\in{\mathbb{SY}}_n$, we see that 
$$
X_\mu\supset X_\lambda\ \Longleftrightarrow\ 
Y_\mu\supset Y_\lambda\ \Longleftrightarrow\ 
\mu_j\leq \lambda_j \mbox{ for } 1\leq j\leq n.
$$
Define the {\em Bruhat order} by $\mu\leq \lambda$ if 
$\mu_j\leq \lambda_j$ for $1\leq j\leq n$.
Note that $\mu\leq \lambda$ if and only if $\mu_j\leq \lambda_j$ for those
$j$ with $0<\mu_j$.

\begin{ex}{\rm 
Suppose $n=4$.
Then $X_{3\,2\,\overline{1}\,\overline{4}}$ consists of 
those $H\in B_4$ such that 
$$
 \dim H\cap\Span{e_3,e_4}\geq 1,\ 
 \dim H\cap\Span{e_2,e_3,e_4}\geq 2,\ \mbox{and}
 \dim H\cap\Span{e_{\overline{1}},\ldots,e_4}\geq 3.
$$
}
\end{ex}

Define $P_\lambda:=[X_\lambda]$, the cohomology class Poincar{\'e} dual to
the fundamental cycle of $X_\lambda$ in the homology of $B_n$.\
Likewise set $Q_\lambda:=[Y_\lambda]$.
Since Schubert varieties are closures of cells from a decomposition
into (real) even-dimensional cells, these {\em Schubert classes} 
$\{P_\lambda\}$, $\{Q_\lambda\}$ form bases for integral cohomology:
$$
H^*B_n\ =\ \bigoplus_{\lambda} P_\lambda \cdot {\mathbb Z}
\qquad \mbox{and}\qquad
H^*C_n\ =\ \bigoplus_{\lambda} Q_\lambda \cdot {\mathbb Z}.
$$

Each $\lambda\in{\mathbb{SY}}_n$ determines and is determined by its 
diagram, also denoted $\lambda$.
The diagram of $\lambda$ is a left-justified array of $|\lambda|$ boxes 
with $\lambda_j$ boxes in the $j$th row, for $\lambda_j>0$.
Thus
$$
3\,2\,\overline{1}\,\overline{4} \ \longleftrightarrow\ 
\setlength{\unitlength}{.9pt}\begin{picture}(30,20)(0,7)
\put( 0, 0){\line(0,1){20}}\put(0, 0){\line(1,0){20}}
\put(10, 0){\line(0,1){20}}\put(0,10){\line(1,0){30}}
\put(20, 0){\line(0,1){20}}\put(0,20){\line(1,0){30}}
\put(30,10){\line(0,1){10}}\end{picture}\, 
\qquad\mbox{ and }\qquad
4\,2\,1\,\overline{3} \ \longleftrightarrow\ 
\setlength{\unitlength}{.9pt}\begin{picture}(40,20)(0,12)
\put( 0, 0){\line(0,1){30}}\put(0,30){\line(1,0){40}}
\put(10, 0){\line(0,1){30}}\put(0,20){\line(1,0){40}}
\put(20,10){\line(0,1){20}}\put(0,10){\line(1,0){20}}
\put(30,20){\line(0,1){10}}\put(0, 0){\line(1,0){10}}
\put(40,20){\line(0,1){10}}\end{picture}\, .
\raisebox{-15pt}{\rule{0pt}{5pt}}
$$
The Bruhat order corresponds to inclusion of diagrams.
Given $\mu\leq \lambda$, let $\lambda/\mu$ be their set-theoretic
difference.
For instance,
$$
4\,2\,1\,\overline{3}/3\,2\,\overline{1}\,\overline{4}\ \longleftrightarrow\ 
\setlength{\unitlength}{.9pt}\begin{picture}(40,20)(0,12)
\put(30,20){\line(0,1){10}}\put(30,30){\line(1,0){10}}
\put(40,20){\line(0,1){10}}\put(30,20){\line(1,0){10}}
\put( 0, 0){\line(0,1){10}}\put( 0,10){\line(1,0){10}}
\put(10, 0){\line(0,1){10}}\put( 0, 0){\line(1,0){10}}
\put( 0,20){\dashbox{2}(10,10)[t]{}}\put( 0,10){\dashbox{2}(10,10)[t]{}}
\put(10,20){\dashbox{2}(10,10)[t]{}}\put(10,10){\dashbox{2}(10,10)[t]{}}
\put(20,20){\dashbox{2}(10,10)[t]{}}
\end{picture}
\qquad\mbox{ and }\qquad
3\,2\,\overline{1}\,\overline{4}/
1\,\overline{2}\,\overline{3}\,\overline{4}\ \longleftrightarrow\ 
\setlength{\unitlength}{.9pt}\begin{picture}(30,20)(0,7)
\put( 0, 0){\line(0,1){10}}\put( 0, 0){\line(1,0){20}}
\put(10, 0){\line(0,1){20}}\put( 0,10){\line(1,0){30}}
\put(20, 0){\line(0,1){20}}\put(10,20){\line(1,0){20}}
\put(30,10){\line(0,1){10}}\put( 0,10){\dashbox{2}(10,10)[t]{}}
\end{picture}\, .
\raisebox{-15pt}{\rule{0pt}{5pt}}
$$
Two boxes are connected if they share a vertex or an edge; this defines 
{\em components} of $\lambda/\mu$.
We say $\lambda/\mu$ is a {\em skew row} if
$\lambda_1\geq\mu_1\geq\lambda_2\geq\cdots\geq\mu_n$ equivalently, 
if  $\lambda/\mu$ has at most one box in each column.
Thus $4\,2\,1\,\overline{3}/3\,2\,\overline{1}\,\overline{4}$
is a skew row, but 
$3\,2\,\overline{1}\,\overline{4}/
1\,\overline{2}\,\overline{3}\,\overline{4}$
is not.

The {\em special Schubert class} $p_m\in H^*B_n$ ($q_m\in H^*C_n$) is the
class whose diagram consists of a single row of length $m$.
Hence, $p_2 = P_{2\,\overline{1}\,\overline{3}\,\overline{4}}$.
A {\em special Schubert variety} $X_K$ ($Y_K$) is the collection of all
maximal isotropic subspaces which meet a fixed isotropic subspace $K$
nontrivially.
If $\dim K=n+1-m$, then $[X_K]=p_m$ and $[Y_K]=q_m$.

\begin{thm}[Pieri-type Formula]\label{thm:pieri}
For any $\mu\in{\mathbb{SY}}_n$ and $1\leq m\leq n$,
\begin{enumerate}
\item ${\displaystyle
P_\mu\cdot p_m\ =\ \sum_{\lambda/\mu\ \mbox{\scriptsize skew row}}
2^{\delta(\lambda/\mu)-1}\, P_\lambda}$ \qquad  and
\item ${\displaystyle
Q_\mu\cdot q_m\ =\ \sum_{\lambda/\mu\ \mbox{\scriptsize skew row}}
2^{\varepsilon(\lambda/\mu)}\, Q_\lambda}$,
\end{enumerate}
where $\delta(\lambda/\mu)$ counts the components of
the diagram $\lambda/\mu$ and 
$\varepsilon(\lambda/\mu)$ counts the components of $\lambda/\mu$
which do not contain a box in the first column.
\end{thm}

\begin{ex}\label{ex:prod}
{\em 
For example,
\begin{eqnarray*}
P_{3\,2\,\overline{1}\,\overline{4}}\cdot p_2 &=&
2\cdot P_{4\,2\,1\,\overline{3}}\ \,+
\ \ \, P_{4\,3\,\overline{2}\,\overline{1}}\qquad\mbox{ and}\\
Q_{3\,2\,\overline{1}\,\overline{4}}\cdot q_2 &=&
2\cdot Q_{4\,2\,1\,\overline{3}}\ +\ 2\cdot 
Q_{4\,3\,\overline{2}\,\overline{1}}.
\end{eqnarray*}
}
\end{ex}

For $\lambda,\mu,\nu\in{\mathbb{SY}}_n$, there exist integral
constants  
$g^\lambda_{\mu,\nu}$ and $h^\lambda_{\mu,\nu}$ defined by the identities
$$
P_\mu\cdot P_\nu\ =\ \sum_\lambda g^\lambda_{\mu,\nu}\,P_\lambda
\qquad \mbox{and}\qquad
Q_\mu\cdot Q_\nu\ =\ \sum_\lambda h^\lambda_{\mu,\nu}\,Q_\lambda.
$$
These constants were first given a combinatorial formula by
Stembridge~\cite{Stembridge_shifted}.\smallskip

Define $\lambda^c$ by 
$\lambda^c_j:=\overline{\lambda_{n+1-j}}$.
Let $[\mbox{pt}]$ be the class dual to a
point.
The Schubert basis is self-dual with respect to the intersection pairing:
If $|\lambda|=|\mu|$, then 
\begin{equation}\label{eq:poincare}
P_\mu\cdot P_{\lambda^c}
\ =\ 
Q_\mu\cdot Q_{\lambda^c}
\ =\ \left\{\begin{array}{ll}
[\mbox{pt}]&\ \ \mbox{if}\ \lambda=\mu\\
      0    &\ \ \mbox{otherwise}\end{array}\right..
\end{equation}

Define the Schubert variety $X'_{\lambda^c}$ to be 
$$
\{H\in B_n\mid 
\dim H\cap\Span{e_{\overline{n}},\ldots,e_{\lambda_j}}\geq n+1-j
\ \mbox{for}\ 1\leq j \leq n\}.
$$
This is a translate of $X_{\lambda^c}$ by an element of 
$\mbox{\em So}_{2n+1}{\mathbb C}$.
In a similar fashion, define $Y'_{\lambda^c}$, a translate of
$Y_{\lambda^c}$ by an element of $\mbox{\em Sp}_{2n}{\mathbb C}$.
For any $\lambda,\mu$, $X_\mu\bigcap X'_{\lambda^c}$ is generically
transverse~\cite{Kleiman}. 
This is because if $X_\mu$ and $X'_{\lambda^c}$ are {\em any} Schubert
varieties in general position, then there is a basis for $V$
such that these varieties and the quadratic form $\beta$ are as given.

We see that to establish the Pieri-type formula, it suffices to compute the
degree of 
$$
X_\mu\bigcap X'_{\lambda^c}\bigcap X_K \quad\mbox{and}\quad
Y_\mu\bigcap Y'_{\lambda^c}\bigcap Y_K
$$
where $K$ is a general isotropic ($n+1-m$)-plane and 
$|\lambda|=|\mu|+m$.\smallskip

We only do the (more difficult) orthogonal case of Theorem~\ref{thm:pieri}
in full, and indicate the differences for the symplectic case.
We first determine when 
$X_\mu\bigcap X'_{\lambda^c}\neq\emptyset$.
Let $\mu,\lambda\in{\mathbb{SY}}_n$.
Then, by the definition of Schubert varieties, 
$H\in X_\mu\bigcap X'_{\lambda^c}$ implies
$\dim H\bigcap \Span{e_{\mu_j},\ldots,e_{\lambda_j}}\geq 1$,
for every $1\leq j\leq n$.
Hence $\mu\leq\lambda$ is necessary for $X_\mu\bigcap X'_{\lambda^c}$
to be nonempty.
In fact,
$$
X_\mu \bigcap X'_{\lambda^c}\ =\ \left\{
\begin{array}{ll}
\Span{e_{\lambda_1},\ldots,e_{\lambda_n}}&\mbox{ if } \lambda=\mu\\
\emptyset &\mbox{ otherwise},
\end{array}\right.
$$
which establishes~(\ref{eq:poincare}).

Suppose $\mu\leq \lambda$ in ${\mathbb{SY}}_n$.
For each component $d$ of $\lambda/\mu$, let col$(d)$ index the 
columns of $d$ together with the column just to the left of $d$,
which is $0$ if $d$ meets the first column, in that it has a box 
in the first column.
For each component $d$ of $\lambda/\mu$, define a quadratic form $\beta_d$:
$$
\beta_d\ :=\ 
\sum_{j\in \mbox{\scriptsize col}(d)}x_j x_{\overline{\jmath}},
$$
where $x_{\overline{n}},\ldots,x_n$ are coordinates for $V$ dual to the 
basis $e_{\overline{n}},\ldots,e_n$.
For each {\em fixed point} of $\lambda/\mu$ 
($j$ such that $\lambda_j=\mu_j$), define the linear form
$\alpha_j:=x_{\overline{\lambda_j}}$.
If there is no component meeting the first
column, then we say that $0$ is a fixed point of $\lambda/\mu$ and 
define $\alpha_0:=x_0$.
Let $Z_{\lambda/\mu}$ be the common zero locus of these forms $\alpha_j$ and
$\beta_d$.

\begin{lemma}\label{lemma:vanish}
Let $\mu\leq\lambda$ and $H\in X_\mu\bigcap X'_{\lambda^c}$.
Then $H\subset Z_{\lambda/\mu}$.
\end{lemma}

Let ${\mathcal Q}$ be the isotropic points in $V$, the
zero locus of $\beta$.
For each $0\leq i\leq n$, 
there is a unique form among the $\alpha_j$, $\beta_d$ in which one (or
both) of the coordinates $x_{i},x_{\overline{\imath}}$ appears.
Thus $\beta$ is in the ideal generated by these forms $\alpha_j$, 
$\beta_d$ and we see that 
they are dependent on ${\mathcal Q}$.
However, if $\delta=\delta(\lambda/\mu)$ counts the components of
$\lambda/\mu$ and $\varphi$ the number of fixed points, then the collection
of $\varphi$ forms $\alpha_j$   and $\delta-1$ of
the forms $\beta_d$ {\em are} independent on ${\mathcal Q}$.
Moreover, Lemma~\ref{lem:columns} shows that 
$$
n +1 \ =\ \varphi+\delta+\#\mbox{\rm columns of $\lambda/\mu$}.
$$
Thus, if $m=|\lambda|-|\mu|$, then 
$\varphi+\delta-1\leq n-m$, with equality only when $\lambda/\mu$ is a
skew row.
Since ${\mathcal Q}$ has dimension $2n$,
it follows that a general isotropic ($n+1-m$)-plane $K$
meets $Z_{\lambda/\mu}$ only if
$\lambda/\mu$ is a skew row.
We deduce

\begin{thm}\label{thm:geom}
Let $\mu,\lambda\in{\mathbb{SY}}_n$ and suppose $K$ is a general isotropic
$(n+1-m)$-plane with $|\mu|+m = |\lambda|$.
Then
$$
X_\mu \bigcap X'_{\lambda^c}\bigcap X_K\ \neq \ \emptyset
$$
only if $\mu\leq \lambda$ and $\lambda/\mu$ is a skew row.
\end{thm}

Under the hypotheses of Theorem~\ref{thm:geom}, the forms $\alpha_j$ and
$\beta_d$ define $2^{\delta(\lambda/\mu)-1}$ isotropic lines in $K$.
Theorem~\ref{thm:unique} asserts that a general isotropic line in 
$Z_{\lambda/\mu}$ determines a unique
$H\in X_\mu \bigcap X'_{\lambda^c}$, which completes the proof of
Theorem~\ref{thm:pieri}. 

\begin{ex}
{\em 
We show that if  $K\subset {\mathcal Q}$ is a general 3-plane, then 
$$
\#\, X_{3\,2\,\overline{1}\,\overline{4}}\bigcap
X'_{(4\,2\,1\,\overline{3})^c} \bigcap X_K\ =\ 2.
$$
Note that 2 is the coefficient of $P_{4\,2\,1\,\overline{3}}$ in the product
$P_{3\,2\,\overline{1}\,\overline{4}}\cdot p_2$ of Example~\ref{ex:prod}.

First, the local coordinates for 
$X_{3\,2\,\overline{1}\,\overline{4}}\bigcap
X'_{(4\,2\,1\,\overline{3})^c}$
described in Lemma~\ref{lem:loc_coords} show that, for any 
$x,z\in{\mathbb{C}}$, the row span $H$ of the matrix
with rows $g_i$ and columns $e_j$
$$
\begin{array}{l|cccc|c|cccc}
{} &e_{\overline{4}}&e_{\overline{3}}
&e_{\overline{2}}&e_{\overline{1}}&e_0&e_1&e_2&e_3&e_4\\
\hline
g_1&0&0&0&0&0&0&0&-x&1\\
g_2&0&0&0&0&0&0&1&0&0\\
g_3&0&0&0&1&2z&-2z^2&0&0&0\\
g_4&x&1&0&0&0&0&0&0&0
\end{array}
$$
is in $X_{3\,2\,\overline{1}\,\overline{4}}\bigcap
X'_{(4\,2\,1\,\overline{3})^c}$.
Suppose $K$ is the row span of the matrix with rows $v_i$
$$
\begin{array}{l|cccc|c|cccc}
{} &e_{\overline{4}}&e_{\overline{3}}
&e_{\overline{2}}&e_{\overline{1}}&e_0&e_1&e_2&e_3&e_4\\
\hline
v_1&0&1&0&1&0& 0&1&0&1\\
v_2&1&1&0&1&2&-2&1&1&-1\\
v_3&0&0&1&0&0&-1&0&0&0
\end{array}
$$
Then $K$ is an isotropic 3-plane, and the forms
\begin{eqnarray*}
\beta_0&=& 2x_{\overline{1}}x_1 + x_0^2\\
\beta_d&=&x_{\overline{4}}x_4 + x_{\overline{3}}x_3\\
\alpha_2&=& x_{\overline{2}}
\end{eqnarray*}
define the 2 isotropic lines $\Span{v_1}$ and $\Span{v_2}$ in $K$.
Lastly, for $i=1,2$, there is a unique 
$H_i\in  X_{3\,2\,\overline{1}\,\overline{4}}\bigcap
X'_{(4\,2\,1\,\overline{3})^c}$
with $v_i\in H_i$.
In these coordinates,
$$
H_1\ :\  x=z=0\quad\mbox{and}\quad
H_2\ :\ x=z=1.
$$
}
\end{ex}

In  the symplectic case,  isotropic $K$ are not
contained in a quadric ${\mathcal Q}$, the form $\alpha_0=x_0$ does not
arise,  only components which do not meet the first column give
quadratic forms $\beta_d$, 
and  the analysis of Lemma~\ref{lem:components}~(2) in
Section~\ref{sec:triple} is (slightly) different.

\section{The intersection of two Schubert varieties}

We study the intersection $X_\mu\bigcap X'_{\lambda^c}$ of two Schubert
varieties.
Our main result, Theorem~\ref{thm:main}, expresses 
$X_\mu\bigcap X'_{\lambda^c}$ as a product whose factors correspond to
components of $\lambda/\mu$, and each factor is itself an intersection of
two Schubert varieties.
These factors are described in Lemmas~\ref{lemma:0comp} and~\ref{lem:offdiag},
and in Corollary~\ref{cor:classical}.
These are crucial to the proof of the Pieri-type formula that we complete in
Section 3.
Also needed is Lemma~\ref{lem:subspace}, which identifies a particular
subspace of $H\cap \Span{e_1,\ldots,e_n}$ for 
$H\in X_\mu\bigcap X'_{\lambda^c}$.

For Lemma~\ref{lem:subspace}, we work in the
(classical) Grassmannian $G_k(V^+)$ of $k$-planes in 
$V^+:=\Span{e_1,\ldots,e_n}$.
For basic definitions and results see any of
~\cite{Hodge_Pedoe,Griffiths_Harris,Fulton_tableaux}. 
Schubert subvarieties $\Omega_\sigma,\Omega'_{\sigma^c}$ of 
$G_k(V^+)$  are indexed by partitions
$\sigma\in{\mathbb{Y}}_k$, where
$n - k\geq \sigma_1\geq\cdots\geq\sigma_k\geq0$.
For $\sigma\in{\mathbb{Y}}_k$ define $\sigma^c\in{\mathbb{Y}}_k$ by
$\sigma^c_j=n-k-\sigma_{k+1-j}$.
For $\sigma,\tau\in{\mathbb{Y}}_k$, define
\begin{eqnarray*}
\Omega_\tau&:=& \{H\in G_k(V^+)\mid
 \dim H\cap\Span{e_{k+1-j+\tau_j},\ldots,e_n}\geq j,\ 1\leq j\leq k\}\\
\Omega'_{\sigma^c}&:=& \{H\in G_k(V^+)\mid
 \dim H\cap\Span{e_1,\ldots,e_{j+\sigma_{k+1-j}}}\geq j,\ 1\leq j\leq k\}.
\end{eqnarray*}

Let $\lambda,\mu\in\mathbb{SY}_n$ with $\mu\leq \lambda$, and 
suppose $\mu_k>0>\mu_{k+1}$.
Define partitions $\sigma$ and $\tau$ in ${\mathbb{Y}}_k$
(which depend upon $\lambda$ and $\mu$)  by
\begin{eqnarray*}
\tau&:=& \mu_1-k\geq \cdots\geq \mu_k-1\geq 0\\
\sigma&:=& \lambda_1-k\geq \cdots\geq \lambda_k-1>0.
\end{eqnarray*}

\begin{lemma}\label{lem:subspace}
Let $\mu\leq \lambda\in{\mathbb{SY}}_n$, and define 
$\sigma,\tau\in{\mathbb{Y}}_k$, and $k$ as above.
If $H\in X_\mu\bigcap X'_{\lambda^c}$, then 
$H\cap V^+$ contains a $k$-plane 
$L\in\Omega_\tau\bigcap\Omega'_{\sigma^c}$.
\end{lemma}

\noindent{\bf Proof. }
Suppose first that $H\in X_\mu$ with 
$\dim H\cap \Span{e_{\mu_j},\ldots,e_n}=j$ for $j=k$ and $k+1$.
Since $\mu_k>0>\mu_{k+1}$, we see that 
$L:= H\cap V^+$ has dimension $k$.
If $H\in X'_{\lambda^c}$ in addition, it
is an exercise in the definitions to verify that 
$L\in\Omega_\tau\bigcap\Omega'_{\sigma^c}$.
The lemma follows as such $H$ are dense in $X_\mu$.
\QED

The first step towards Theorem~\ref{thm:main}
is the following combinatorial lemma.

\begin{lemma}\label{lem:columns}
Let $\varphi$ count the fixed points and $\delta$
the components of $\lambda/\mu$.
Then
$$
n +1\ =\ \varphi+\delta+\#\mbox{\rm columns of $\lambda/\mu$},
$$
and $\mu_j>\lambda_{j+1}$ precisely when $|\mu_j|$ is an empty column of
$\lambda/\mu$.
\end{lemma}

\noindent{\bf Proof. }
Suppose $k$ is a column not meeting $\lambda/\mu$.
Thus, there is no $i$ for which $\mu_i<k\leq\lambda_i$.
Let $j$ be the index such that $|\mu_j|=k$.
If $\mu_j=k$, then we must also have $\lambda_{j+1}<k$, 
as $\mu_{j+1}<k$.
Either $\mu_j=\lambda_j$ is a fixed point of
$\lambda/\mu$ or else $\mu_j<\lambda_j$, so that $k$ is the column
immediately to the left of a component $d$ which does not meet the first
column. 

If $\mu_j=-k$, then $\lambda_j=-k$, for otherwise $k=\lambda_i$ for some
$i$, and for this $i$ we must necessarily have $\mu_i<k$, 
contradicting $k$ being an empty column.
This proves the lemma, as $0$ is either a fixed point of $\lambda/\mu$ or
else $\lambda/\mu$ has a component meeting the first  column, but not
both.\QED

Let $d_0$ be the component of $\lambda/\mu$ meeting the first column (if
any). 
Define mutually orthogonal subspaces 
$V_\varphi,V_0$, and $V_d$, for $d$ a component of $\lambda/\mu$ not
meeting the first column ($0\not\in\mbox{col}(d)$) as follows:
\begin{eqnarray*}
V_\varphi\: &:=& \Span{e_{\mu_j},e_{\overline{\mu_j}}\mid \mu_j=\lambda_j},\\
V_0\:\, &:=& \Span{e_0,e_k,e_{\overline{k}}\mid k\in\mbox{col}(d_0)},\\
V^-_d &:=& \Span{e_k\mid k\in\mbox{col}(d)},\\
V^+_d &:=& \Span{e_{\overline{k}}\mid k\in\mbox{col}(d)},
\end{eqnarray*}
and set $V_d:= V^-_d\oplus V^+_d$.
Then 
$$
V\ =\ V_\varphi\oplus V_0\oplus
\bigoplus_{0\not\in\mbox{\scriptsize col}(d)} V_d.
$$
For each fixed point $\mu_j=\lambda_j$ of $\lambda/\mu$, 
define the linear form 
$\alpha_j:=x_{\overline{\mu_j}}$.
For each component $d$ of $\lambda/\mu$, let the quadratic form $\beta_d$ be
the restriction of the form $\beta$ to $V_d$.
Composing with the projection of $V$ to $V_d$ gives a quadratic form (also
written $\beta_d$) on $V$.
If there is no component meeting the first column,
define $\alpha_0:=x_0$ and call $0$ a fixed point of $\lambda/\mu$.
If $0\not\in\mbox{col}(d)$, then the form $\beta_d$ identifies 
$V^+_d$ and $V^-_d$ as dual vector spaces.

\begin{lemma}\label{lemma:intersection}
Let $H\in X_\mu\bigcap X'_{\lambda^c}$.
Then
\begin{enumerate}
\item $H\bigcap V_\varphi = \Span{e_{\mu_j}\mid\mu_j=\lambda_j}$.
\item $\dim H\bigcap V_0 = \#\mbox{col}(d_0)$.
\item For all components $d$ of $\lambda/\mu$ which do not meet the first
column, 
\begin{eqnarray*}
\dim H\bigcap V^+_d&=& \#\mbox{rows of }d,\\
\dim H\bigcap V^-_d&=& \#\mbox{col}(d) - \#\mbox{rows of }d,
\end{eqnarray*}
and $\left(H\bigcap V^-_d\right)^\perp = H\bigcap V^+_d$.
\end{enumerate}
\end{lemma}

\noindent{\bf Proof of Lemma~\ref{lemma:intersection}. }
Let $H\in X_\mu\bigcap X'_{\lambda^c}$.
Suppose $\mu_j>\lambda_{j+1}$
so that $|\mu_j|$ is an empty column of $\lambda/\mu$.
Then the definitions of Schubert varieties imply
$$
H\ =\ H\cap\Span{e_{\overline{n}},\ldots,e_{\lambda_{j+1}}}\oplus
H\cap\Span{e_{\mu_j},\ldots,e_n}.
$$

Suppose $d$ is a component not meeting the first column.
If the rows of $d$ are $j,\ldots,k$, then 
\begin{eqnarray*}
H\cap V^+_d&=& H\cap\Span{e_{\mu_k},\ldots,e_{\lambda_j}}\\
&=& H  \cap\Span{e_{\overline{n}},\ldots,e_{\lambda_j}}
      \cap\Span{e_{\mu_k},\ldots,e_n},
\end{eqnarray*}
and so has dimension at least $k-j+1$.

Similarly, if $l,\ldots,m$ are the indices $i$ with 
$\overline{\lambda_{j}}\leq \mu_i,\lambda_i\leq \overline{\mu_k}$, 
then 
$H\bigcap V^-_d$ has dimension at least $l-m+1$.
Hence
$\dim V_d/2 = \#\mbox{col}(d)=k+m-l-j+2$, as 
$\lambda_j,\ldots,\lambda_k,\overline{\lambda_l},\ldots,\overline{\lambda_m}$
are the columns of $d$.

Since $H$ is isotropic, $\dim H^+_d + \dim H^-_d \leq \#\mbox{col}(d)$,
which proves the first part of (3).
Moreover, $H\bigcap V^+_d=\left(H\bigcap V^-_d\right)^\perp$:
Since $H$ is isotropic, we have $\subset$, and equality follows by
dimension counting.

Similar arguments prove the other statements.
\QED

For $H\in X_\mu\bigcap X'_{\lambda^c}$, define 
$H_\varphi:= H\bigcap V_\varphi$,
$H_0:= H\bigcap V_0$,
$H^+_d := H\bigcap V^+_d$, and 
$H^-_d := H\bigcap V^-_d$.
Then $H_\varphi\subset V_\varphi$ is the zero locus of the linear 
forms $\alpha_j$, $H_0$ is isotropic in $V_0$, and, for each component $d$
of $\lambda/\mu$ not meeting the first column, 
$H_d:= H^+_d\oplus H^-_d$ is isotropic in $V_d$, which proves
Lemma~\ref{lemma:vanish}. 
Moreover, $H$ is the orthogonal direct sum of $H_\varphi$, $H_0$, and the 
$H_d$.

\begin{thm}\label{thm:main}
The map 
$$
\{H_0\mid H\in X_\mu\bigcap X'_{\lambda^c}\}\  \times
\prod_{0\not\in\mbox{\scriptsize col}(d)}
\{H_d\mid H\in X_\mu\bigcap X'_{\lambda^c}\}
\longrightarrow\  X_\mu\bigcap X'_{\lambda^c}
$$
defined by
$$
(H_0,\ldots,H_d,\ldots) \longmapsto 
\Span{H_\varphi,H_0,\ldots,H_d,\ldots}
$$
is an isomorphism.
\end{thm}

\noindent{\bf Proof. }
By the previous discussion, it is an injection.
For surjectivity, note that both sides have the same dimension.
Indeed, $\dim  X_\mu\bigcap X'_{\lambda^c}=|\lambda|-|\mu|$,
the number of boxes in $\lambda/\mu$.
Lemmas~\ref{lemma:0comp} and \ref{lem:offdiag}
show that the factors of the domain each have dimension equal to the
number of boxes in the corresponding components.
\QED

Suppose there is a component, $d_0$, meeting the first column.
Let $l$ be the largest column in $d_0$, and
define $\lambda(0),\mu(0)\in{\mathbb{SY}}_l$ as follows:
Let $j$ be the first row of $d_0$ so that $l=\lambda_j$.
Then, since $d_0$ is a component, for each $j\leq i<j+l-1$, we have 
$\lambda_{i+1}\geq \mu_i$ and $l=\overline{\mu_{j+l-1}}$.
Set 
\begin{eqnarray*}
\mu(0)&:=& \mu_j>\cdots >\mu_{j+l-1}\\
\lambda(0)&:=& \lambda_j>\cdots >\lambda_{j+l-1}
\end{eqnarray*}
Define $\lambda(0)^c$ by 
$\lambda(0)^c_p:= \overline{\lambda(0)_{l+1-p}} 
= \overline{\lambda_{j+l-p}}$.
The following lemma is a straightforward consequence of these definitions.

\begin{lemma}\label{lemma:0comp}
With the above definitions,
$$
\{H_0\mid H\in X_\mu\bigcap X'_{\lambda^c}\}\ =\ 
X_{\mu(0)}\bigcap X'_{\lambda(0)^c},
$$
as subvarieties of $B_k\simeq B(V_0)$, and $\lambda(0)/\mu(0)$ has a unique
component meeting the first column and no fixed points. 
\end{lemma}

We similarly identify 
$\{H_d\mid H\in X_\mu\bigcap X'_{\lambda^c}\}$
as an intersection $X_{\mu(d)}\cap X'_{\lambda(d)^c}$ of Schubert varieties
in $B_{\#\mbox{\scriptsize col}(d)}\simeq B(\Span{e_0,V_d})$.
Let  $j,\ldots,k$ be the rows of $d$ and $l,\ldots,m$ be the indices $i$ 
with  $\overline{\lambda_j}\leq \mu_i,\lambda_i\leq \overline{\mu_k}$, 
as in the
proof of Lemma~\ref{lemma:intersection}.
Let $p=\#\mbox{col}(d)$ and define
$\lambda(d),\mu(d)\in {\mathbb{SY}}_p$ as follows.
Set $a = \mu_k$, and define 
\begin{eqnarray*}
\mu(d)&:=& \mu_j-a+1>\cdots>
\hspace{24pt}1\hspace{24pt}
>\mu_l+a-1>\cdots>\mu_m+a-1\\
\lambda(d)&:=& \lambda_j-a+1>\cdots>\lambda_k-a+1>
\lambda_l+a-1>\cdots>\lambda_m+a-1
\end{eqnarray*}
As with Lemma~\ref{lemma:0comp}, the following lemma is straightforward.

\begin{lemma}\label{lem:offdiag}
With these definitions,
$$
\{H_d\mid H\in X_\mu\bigcap X'_{\lambda^c}\}\ \simeq\ 
X_{\mu(d)}\bigcap X'_{\lambda(d)^c}
$$
as subvarieties of $B_p\simeq B(\Span{e_0,V_d})$
and $\lambda(d)/\mu(d)$ has a unique component not meeting the first column
and no fixed points.
\end{lemma}

Suppose now that $\mu,\lambda\in{\mathbb{SY}}_n$ where $\lambda/\mu$ has a
unique component $d$ not meeting the first column and no fixed points.
Suppose $\lambda$ has $k$ rows.
A consequence of Lemma~\ref{lemma:intersection} is that the map 
$H^+_d \mapsto \Span{H^+_d,\left(H^+_d\right)^\perp}$ gives an isomorphism
\begin{equation}\label{isomorphism}
	\{ H^+_d\mid H\in X_\mu\bigcap X'_{\lambda^c}\}\ 
	\stackrel{\sim}{\longrightarrow}\ 
	X_\mu\bigcap X'_{\lambda^c}.
\end{equation}

The following corollary of Lemma~\ref{lem:subspace} identifies the domain.

\begin{cor}\label{cor:classical}
With $\mu,\lambda$ as above and $\sigma,\tau$, and $k$ as defined in the 
paragraph
preceding Lemma~\ref{lem:subspace}, we have:
$$
\{H^+_d\mid H\in X_\mu\bigcap X'_{\lambda^c}\}\ =\ 
\Omega_{\tau}\bigcap\Omega_{\sigma^c},
$$
as subvarieties of $G_k(V^+)$.
\end{cor}

\begin{rem}\label{rem:similarities}
{\em 
The symplectic analogs of Lemma~\ref{lem:offdiag} and
Corollary~\ref{cor:classical}, which are identical save for the necessary
replacement of $Y$ for $X$ and $C_p$ for $B_p$, show an unexpected
connection between the geometry of the  
symplectic and orthogonal Grassmannians.
Namely, suppose $\lambda/\mu$ has no component meeting the diagonal.
Then the projection map
$V \twoheadrightarrow W$ defined by 
$$
e_i\ \longmapsto\ \left\{\begin{array}{ll}0&\ \mbox{if}\  i=0\\
f_i&\ \mbox{otherwise}
\end{array}\right.
$$
and its right inverse $W\hookrightarrow V$ defined by
$f_j\mapsto e_j$ induce isomorphisms
$$
X_\mu\bigcap X'_{\lambda^c}\ \stackrel{\sim}{\longleftrightarrow}\ 
Y_\mu\bigcap Y'_{\lambda^c}.
$$
}
\end{rem}

\section{Pieri-type intersections of Schubert varieties}\label{sec:triple}

Fix $\lambda/\mu$ to be a skew row with $|\lambda|-|\mu|=m$.
Let $Z_{\lambda/\mu}\subset{\mathcal Q}$ be the zero locus of the 
forms $\alpha_j$ and $\beta_d$ of \S 2.
If $\lambda/\mu$ has $\delta$ components, then 
as a subvariety of ${\mathcal Q}$,  $Z_{\lambda/\mu}$ is 
the generically transverse intersection of the zero loci of the forms
$\alpha_j$ and any $\delta-1$ of the forms $\beta_d$.
It follows that a general $(n+1-m)$-plane $K\subset{\mathcal Q}$ meets
$Z_{\lambda/\mu}$ in $2^{\delta-1}$ lines. 
Thus if $\Span{v}\subset Z_{\lambda/\mu}$ is a general line, then
$$
\# X_\mu\bigcap X'_{\lambda^c}\bigcap X_K \ =\ 
2^\delta\cdot \# X_\mu\bigcap X'_{\lambda^c}\bigcap X_{\Span{v}}.
$$
Theorem~\ref{thm:pieri} is a consequence of this observation and
the following:

\begin{thm}\label{thm:unique}
Let $\lambda/\mu$ be a skew row, $Z_{\lambda/\mu}$ be as above, 
and $\Span{v}$ a general
line in $Z_{\lambda/\mu}$.
Then $X_\mu\bigcap X'_{\lambda^c}\bigcap X_{\Span{v}}$ is a singleton.
\end{thm}

\noindent{\bf Proof. }
Let ${\mathcal Q}_0$ be the cone of isotropic points in $V_0$ and 
${\mathcal Q}_d$ the cone of isotropic points in $V_d$ for 
$d\neq d_0$.
Since 
$$
Z_{\lambda/\mu}\ =\ H_\varphi \oplus {\mathcal Q}_0\oplus
\bigoplus_{0\not\in\mbox{\scriptsize col}(d)} {\mathcal Q}_d,
$$
we see that a general non-zero vector $v$ in $Z_{\lambda/\mu}$ has the form
$$
v\ =\ \sum_{\mu_j=\lambda_j} a_je_{\mu_j}\ +\ v_0\ +
\sum_{0\not\in\mbox{\scriptsize col}(d)} v_d,
$$
where $a_j\in {\mathbb C}^\times$ and $v_0\in {\mathcal Q}_0$,
$v_d\in {\mathcal Q}_d$ are general vectors.

Thus, if $H\in X_\mu\bigcap X'_{\lambda^c}\bigcap X_{\Span{v}}$, we see that
$v_0\in H_0$ and $v_d\in H_d$.
By Theorem~\ref{thm:main},  $H$ is determined by $H_0$ and the $H_d$,
thus it suffices to prove that $H_0$ and the $H_d$ are uniquely determined.
The identifications of Lemmas~\ref{lemma:0comp} and~\ref{lem:offdiag} show
that this is just the case of the theorem when $\lambda/\mu$ has a single
component, which is Lemma~\ref{lem:components} below.
\QED

\begin{lemma}\label{lem:components}
Suppose $\lambda,\mu\in{\mathbb{SY}}_n$ where $\lambda/\mu$ is a skew row
with a unique component and no fixed points.
\begin{enumerate}
\item
If $\lambda/\mu$ does not meet the first column and $v\in {\mathcal Q}_d$ is
a general vector, then $X_\mu\bigcap X'_{\lambda^c}\bigcap X_{\Span{v}}$
is a singleton.
\item
If $\lambda/\mu$ meets the first column and $v\in {\mathcal Q}$ is
general, then $X_\mu\bigcap X'_{\lambda^c}\bigcap X_{\Span{v}}$
is a singleton.
\end{enumerate}
\end{lemma}

\noindent{\bf Proof  of (1). }
Let $v\in{\mathcal Q}_d$ be a general vector.
Since ${\mathcal Q}_d\subset V^+\oplus V^-$, 
$v=v^+\oplus v^-$ with $v^+ \in V^+$ and $v^-\in V^-$.
Suppose $\mu_k>0>\mu_{k+1}$.
Consider the set
$$
\{H^+\in G_k(V^+)\mid v\in H^+\oplus\left(H^+\right)^\perp\}
\ =\ \{H^+\mid v^+\in H^+\subset (v^-)^\perp\}.
$$
This is a Schubert variety 
$\Omega''_{h(n-k,k)}$ of $G_kV^+$, where $h(n-k,k)$ is the partition of hook
shape with a single row of length $n-k$ and a single 
column of length $k$.

Under the isomorphisms of (\ref{isomorphism}) and Lemma~\ref{lem:offdiag},
and with the identification  
of Corollary~\ref{cor:classical}, we see that 
$$
X_\mu\bigcap X'_{\lambda^c}\bigcap X_{\Span{v}}\ \simeq\ 
\Omega_\tau\bigcap\Omega'_{\sigma^c}\bigcap\Omega''_{h(n-k,k)},
$$
where $\sigma,\tau$ are as defined in the paragraph preceding 
Corollary~\ref{cor:classical}.
For $\rho\in{\mathbb{Y}}_k$, let $S_\rho:=[\Omega_\rho]$ be the cohomology
class Poincar\'e dual to the fundamental cycle of $\Omega_\rho$ in
$H^*G_kV^+$.
The multiplicity we wish to compute is 
\begin{equation}\label{Schur:calc}
	\deg (S_\tau\cdot S_{\sigma^c}\cdot S_{h(n-k,k)}).
\end{equation}
By a double application of the classical Pieri's formula 
(as $S_{h(n-k,k)}=S_{n-k}\cdot S_{1^{k-1}}$), we see that 
(\ref{Schur:calc}) is either 1 or 0, depending upon whether or not
$\sigma/\tau$ has exactly one box in each diagonal.
But this is the case, as the transformation
$\mu,\lambda \longmapsto \tau,\sigma$
takes columns to diagonals.
\QED

Our proof of Lemma~\ref{lem:components} (2) uses a system of 
local coordinates for 
$X_\mu\bigcap X'_{\lambda^c}$.
Let $\lambda/\mu$ be as in Lemma~\ref{lem:components} (2), and suppose
$\lambda_{k+1}=1$.
For $y_2,\ldots,y_n,x_0,\ldots,x_{n-1}\in{\mathbb C}$, 
define vectors $g_j\in V$ as follows:
\begin{equation}\label{loc_coords}
g_j\ := \ \left\{\begin{array}{ll}
{\displaystyle 
e_{\lambda_j} + \sum_{i=\mu_j}^{\lambda_j-1}x_i\,e_i}&\ j\leq k\\
{\displaystyle 
-2x_0^2e_1 + 2x_0e_0+ e_{\overline{1}} 
+ \sum_{i=\mu_{k+1}}^{\overline{2}}y_i\,e_i}&\  j= k+1\\
{\displaystyle 
e_{\lambda_j} + \sum_{i=\mu_j}^{\lambda_j-1}y_i\,e_i}&\  j>k+1\\
\end{array}\right..
\end{equation}

\begin{lemma}\label{lem:loc_coords}
Let $\lambda,\mu\in{\mathbb{SY}}_n$ where $\lambda/\mu$ is a skew row
meeting the first column with no fixed points, and define
$\tau,\sigma\in{\mathbb{Y}}_k$, and $k$ as for Lemma~\ref{lem:subspace}.
Then
\begin{enumerate}
\item For any $x_1,\ldots,x_{n-1}\in{\mathbb C}$, we have
$\Span{g_1,\ldots,g_k}\ \in\ \Omega_\tau\bigcap\Omega'_{\sigma^c}$.
\item
For and $x_0,\ldots,x_{n-1}\in{\mathbb C}$ with 
$x_{\overline{\mu_{k+1}}},\ldots,x_{\overline{\mu_{n-1}}}\neq 0$,
the condition that $H:= \Span{g_1,\ldots,g_n}$ is isotropic
determines a unique 
$H\in X_\mu\bigcap X'_{\lambda^c}$.
\end{enumerate}
Moreover, these coordinates parameterize dense subsets of the
intersections. 
\end{lemma}

\noindent{\bf Proof. }
The first statement is immediate from the definitions.

For the second, note that each $g_j\in{\mathcal Q}$.
The conditions that $\Span{g_1,\ldots,g_n}$ is isotropic are
$$
\beta(g_i,g_j)\ =\ 0\quad \mbox{for}\quad i\leq k< j.
$$

Only $n-1$ of these are not identically zero.
Indeed, for $i\leq k<j$,
$$
\beta(g_i,g_j)\ \not\equiv\  0\ \  \Longleftrightarrow\ \ 
\left\{\begin{array}{ll}
\mbox{either}\ & \overline{\lambda_j}<\mu_i<\overline{\mu_j},\\
\mbox{or}& \mu_i<\overline{\mu_j}<\lambda_i. \end{array}\right.
$$
Moreover, if we order the variables $y_2<\cdots<y_n<x_0<\cdots<x_{n-1}$,
then, in the lexicographic term order, the leading term of $\beta(g_i,g_j)$
for $i\leq k<j$ is
$$
\begin{array}{cl}
y_{\lambda_i}&\ \mbox{if}\ \overline{\lambda_j}<\mu_i<\overline{\mu_j},\\
y_{\overline{\mu_j}}x_{\overline{\mu_j}}
&\ \mbox{if}\ \mu_i<\overline{\mu_j}<\lambda_i,\quad \mbox{or}\\
y_n=y_{\overline{\mu_n}}&\ \mbox{if}\ i=1,\  j=n.\end{array}
$$

Since $\{2,\ldots,n\}=\{\lambda_2,\ldots,\lambda_{k-1},
\overline{\mu_k},\ldots,\overline{\mu_n}\}$,
each $y_l$ appears as the leading term of a unique 
$\beta(g_i,g_j)$ with $i<k\leq j$,
thus these $n-1$ non-trivial equations $\beta(g_i,g_j)=0$ determine
$y_2,\ldots,y_n$ uniquely.

These coordinates parameterize an $n$-dimensional subset of 
$X_\mu\bigcap X'_{\lambda^c}$.
Since $\dim(X_\mu\bigcap X'_{\lambda^c})=n$ and 
$X_\mu\bigcap X'_{\lambda^c}$ is irreducible~\cite{Deodhar}, 
this subset is dense, which completes the proof.
\QED

\begin{ex}
{\em 
Let $\lambda=6\,5\,3\,1\,\overline{2}\,\overline{4}$
and $\mu = 5\,3\,1\,\overline{2}\,\overline{4}\,\overline{6}$
so $k=3$.
We display the components of the vectors $g_i$ in a matrix
$$
\begin{array}{l|cccccc|c|cccccc}
{} &e_{\overline{6}}&e_{\overline{5}}&e_{\overline{4}}&e_{\overline{3}}
&e_{\overline{2}}&e_{\overline{1}}&e_0&e_1&e_2&e_3&e_4&e_5&e_6\\
\hline
g_1&0&0&0&0&0&0&0&0&0&0&0&x_5&1\\
g_2&0&0&0&0&0&0&0&0&0&x_3&x_4&1&0\\
g_3&0&0&0&0&0&0&0&x_1&x_2&1&0&0&0\\
g_4&0&0&0&0&y_2&1&2x_0&-2x_0^2&0&0&0&0&0\\
g_5&0&0&y_4&y_3&1&0&0&0&0&0&0&0&0\\
g_6&y_6&y_5&1&0&0&0&0&0&0&0&0&0&0
\end{array}
$$

Then there are 5 non-zero equations $\beta(g_i,g_j)=0$ with 
$i\leq 3<j$:
\begin{eqnarray*}
0\ =\ \beta(g_3,g_4) &=& y_2 x_2 + x_1\\
0\ =\ \beta(g_3,g_5) &=& y_3+x_2\\
0\ =\ \beta(g_2,g_5) &=& y_4x_4 + y_3x_3\\
0\ =\ \beta(g_2,g_6) &=& y_5+x_4\\
0\ =\ \beta(g_1,g_6) &=& y_6+x_5y_5
\end{eqnarray*}
Solving, we obtain:
$$
y_2=-x_1/x_2,\ 
y_3=-x_2,\ 
y_4=-y_3x_3/x_4,\ 
y_5=-x_4,\ \mbox{and}\ 
y_6=-x_5y_5.
$$
}\end{ex}

\noindent{\bf Proof of Lemma~\ref{lem:components} (2). }
Suppose $\lambda,\mu\in{\mathbb{SY}}_n$ where $\lambda/\mu$ is a skew row
with a single component meeting the first column and no fixed points.
Let $v\in{\mathcal Q}$ be a general vector and consider the condition that
$v\in H$ for $H\in X_\mu\bigcap X'_{\lambda^c}$.
Let $\sigma,\tau\in{\mathbb Y}_k$ be defined as in the paragraph preceding
Lemma~\ref{lem:subspace}. 
We first show that there is a unique 
$L\in\Omega_\tau\bigcap\Omega_{\sigma^c}$ with 
$L\subset H$, and then argue that $H$ is unique.

The conditions on $\mu$ and $\lambda$ imply that 
$\mu_n=\overline{n}$ and $\mu_j=\lambda_{j+1}$ for $j<n$.
We further suppose that $\lambda_{k+1}=1$, so that the last row of 
$\lambda/\mu$ has length 1.
This is no restriction, as the isomorphism of $V$  defined by 
$e_j\mapsto e_{\overline{\jmath}}$ sends 
$X_{\mu}\bigcap X'_{\lambda^c}$ to 
$X_{\lambda^c}\bigcap X'_{(\mu^c)^c}$
and one of $\lambda/\mu$ or $\mu^c/\lambda^c$ has last row of length
1.\smallskip 

Let $v\in {\mathcal Q}$ be general.
If necessary, scale $v$ so that its  $e_{\overline{1}}$-component is
1.
Let $2z$ be its $e_0$-component, then necessarily its
$e_1$-component is $-2z^2$.
Let $v^-\in V^-$ be the 
projection of $v$ to $V^-$.
Similarly define $v^+\in V^+$.
Set
$v':= v^+ + 2z^2e_1$, so that  
$\beta(v^-,v')=0$ and 
$$
v\ =\ v^- + 2z(e_0 - z e_1) + v'.
$$

Let $H\in X_\mu\bigcap X'_{\lambda^c}$, and suppose that $v\in H$.
In the notation of Lemma~\ref{lem:subspace}, let 
$L\in \Omega_\tau\bigcap\Omega_{\sigma^c}$ be a $k$-plane in 
$H\bigcap V^+$.
If $H$ is general, in that 
$$
\dim H\bigcap\Span{e_{\overline{n}},\ldots,e_{\lambda_{k+2}}}\ =\ 
\dim H\bigcap\Span{e_{\overline{n}},\ldots,e_0}\ =\ 
n-k-1,
$$
then  $\Span{L,e_1}$ is the projection of $H$ to $V^+$.
As $v\in H$, we have $v^+\in \Span{L,e_1}$.
Since $L\subset v^\perp\bigcap V^+=(v^-)^\perp$,
we see that $v'\in L$, and hence
$$
v'\ \in\ L\ \subset\  (v^-)^\perp.
$$
As in the proof of part (1), there is a 
(necessarily unique) such $L\in\Omega_\tau\bigcap\Omega_{\sigma^c}$
if and only if $\sigma/\tau$ has a unique box in each diagonal.
But this is the case, as the transformation
$\mu,\lambda \longrightarrow \tau,\sigma$
takes columns (greater than 1) to diagonals.
\smallskip

To complete the proof, we use the local coordinates  for
$X_\mu\bigcap X'_{\lambda^c}$ and $\Omega_\tau\bigcap\Omega_{\sigma^c}$ 
of Lemma~\ref{lem:loc_coords}.
Since $v$ is general, we may assume that that the $k$-plane 
$L\in\Omega_\tau\bigcap\Omega_{\sigma^c}$
determined by $v'\in L\subset(v^-)^\perp$ has non-vanishing
coordinates 
$x_{\overline{\mu_{k+1}}},\ldots,x_{\overline{\mu_{n-1}}}$, 
so that there is an 
$H\in X_\mu\bigcap X'_{\lambda^c}$ in this system of coordinates 
with $L=H\cap V^+$.

Such an $H$ is determined up to a choice of coordinate $x_0$.
The requirement that $v\in H$ forces the projection 
$\Span{e_{\overline{1}}+ 2x_0e_0}$ of $H$ to 
$\Span{e_{\overline{1}},e_0}$ to contain
$e_{\overline{1}}+ 2ze_0$, the projection of $v$ to 
$\Span{e_{\overline{1}},e_0}$.
Hence $x_0=z$, and it follows that there is at most one
$H\in X_\mu\bigcap X'_{\lambda^c}$ with $v\in H$.
Let $g_1,\ldots,g_n$ be the vectors~(\ref{loc_coords}) determined by the
coordinates $x_1,\ldots,x_{n-1}$ for $L$ with $x_0=z$.
We claim $v\in H:=\Span{g_1,\ldots,g_n}$.

Indeed, since $v'\in L$ and 
$v^-\in L^\perp=\Span{g_{k+1}-2z(e_0-ze_1),g_{k+2},\ldots,g_n}$,
we see there exists $\alpha_1,\ldots,\alpha_n\in{\mathbb C}$ with
$$
v^-+v'\ =\ 
\alpha_1 g_1+\cdots + \alpha_{k+1}(g_{k+1}-2z(e_0-ze_1)) 
+\cdots+\alpha_ng_n.
$$
We must have $\alpha_{k+1}=1$, since the $e_{\overline{1}}$-component of
both $v$ and $g_{k+1}$ is 1.
It follows that 
$$
v \ =\ \sum_{i=1}^n \alpha_i g_i\quad \in \ H.\qquad\QED
$$

\noindent{\bf Remarks. }
\begin{enumerate}
\item 
Our desire to give elementary proofs led us to restrict ourselves to the
complex numbers.  
With the appropriate modifications, these same arguments give the same
results for Chow groups of these same varieties over any field of 
characteristic  $\neq 2$.
For example, the appropriate intersection-theoretic constructions and
the properness of a general translate provide a substitute
for our use of transversality.
Then one could argue for the multiplicity of
$2^{\delta-1}$ as follows:

If $\lambda,\mu\in{\mathbb{SY}}_n$ and $K$ is a linear subspace in
${\mathcal{Q}}$, then the scheme-theoretic intersection 
$X_\mu\bigcap X'_{\lambda^c}\bigcap X_K$ is $pr_* \pi^*(K)$, where
$$
\begin{picture}(200,52)
\put(0,0){$K\ \hookrightarrow\ {\mathcal{Q}}$}
\put(80,0){$C_n$}
\put(60,40){$\Xi\ =\ \{(p,H)\mid p\in H\in X_\mu\bigcap X'_{\lambda^c}\}$}
\put(60,35){\vector(-1,-2){11}}
\put(68,35){\vector(1,-2){11}}
\put(46,24){\scriptsize$\pi$}
\put(78,24){\scriptsize$pr$}
\end{picture}
$$
Then intersection theory on the quadric ${\mathcal{Q}}$ (a homogeneous
space) and
Kleiman's Theorem that the intersection with a general translate is
proper~\cite{Kleiman}  
gives a factor of $2^{\delta-1}$ from the intersection multiplicity
of $K$ and the subvariety $Z_{\lambda/\mu}$ of ${\mathcal Q}$ 
consisting of the image of $\pi$.
The arguments of Section 3 show that $\pi$ has degree 1 onto its image.

\item
Conversely, similar to the proof of Lemma~\ref{lem:loc_coords}, we could
give local coordinates for any intersection
$X_\mu\bigcap X'_{\lambda^c}$.
Such a description would enable us to establish transversality directly, and
to dispense with the intersection theory of the classical
Grassmannian.
This would work over any field whose characteristic is not 2,
but would complicate the arguments we gave.

\item We have not investigated to what extent these methods would work
in characteristic 2.

\end{enumerate}

\end{document}